# NONDIFFERENTIABLE FUNCTIONS OF ONE-DIMENSIONAL SEMIMARTINGALES

By George Lowther

We consider decompositions of processes of the form $Y = f(t, X_t)$ where $X$ is a semimartingale. The function $f$ is not required to be differentiable, so Itô's lemma does not apply.

In the case where $f(t,x)$ is independent of $t$, it is shown that requiring $f$ to be locally Lipschitz continuous in $x$ is enough for an Itô-style decomposition to exist. In particular, $Y$ will be a Dirichlet process. We also look at the case where $f(t,x)$ can depend on $t$, possibly discontinuously. It is shown, under some additional mild constraints on $f$, that the same decomposition still holds. Both these results follow as special cases of a more general decomposition which we prove, and which applies to nondifferentiable functions of Dirichlet processes.

Possible applications of these results to the theory of one-dimensional diffusions are briefly discussed.

**1. Introduction.** Suppose that we have a real valued semimartingale $X$ and a function $f : \mathbb{R}_+ \times \mathbb{R} \to \mathbb{R}$. In the case where $f$ is twice continuously differentiable, Itô's lemma shows that $f(t, X_t)$ decomposes as

$$(1) \qquad f(t, X_t) = \int_0^t D_x f(s, X_{s-}) \, dX_s + V_t$$

for a finite variation process $V$. In particular, it follows that $f(t, X_t)$ is itself a semimartingale. The goal of this paper is to generalize this decomposition to situations where $f$ is not differentiable. The case where $f$ is merely once continuously differentiable has been studied previously by several authors and requires going outside the class of semimartingales. Continuous Dirichlet processes were defined by Follmer in [10] as the sum of a continuous local martingale and a process with zero quadratic variation, and it is known that the class of such processes is closed under $C^1$ transformations [2, 4, 7].

These results were applied in [1, 8] and [9] to the study of diffusions with distributional drift.

Noncontinuous Dirichlet processes were defined in [18] as the sum of a semimartingale and a process with zero quadratic variation. It was then shown in [3] that this class of processes is also closed under $C^1$ transformations.

Alternatively, for noncontinuously differentiable functions, decomposition (1) has been studied in [6] assuming that (i) the left derivative $\partial f/\partial t$ exists and is left continuous in $t$ and (ii) there is a decomposition $f = f_h + f_v$ such that $\partial f_h/\partial x$ exists, is continuous and has a left continuous and locally bounded left derivative, and the left derivative $\partial f_v/\partial x$ exists and has a locally bounded variation in $(t, x)$.

In the case where $f(t, x)$ is independent of $t$, we shall show that being locally Lipschitz continuous in $x$ is enough to conclude that the process $V$ in (1) has well-defined quadratic variation with zero continuous part. Working under the slightly generalized definition of a noncontinuous Dirichlet process as the sum of a semimartingale and a process whose quadratic variation has zero continuous parts, this shows that $f(X_t)$ will indeed be a Dirichlet process. We also look at the case where $f$ is a possibly discontinuous function of time. It is required that locally the variation of $f(t, x)$ in $t$ is integrable with respect to $x$. If, additionally, it is locally Lipschitz continuous with left and right derivatives with respect to $x$, then we show that decomposition (1) can be used, and $V$ will have zero continuous quadratic variation. Furthermore, in Section 2 the general situation where $X$ is a Dirichlet process will be looked at. In that case, additional "almost everywhere" differentiability conditions need to be imposed on $f$ and, as we show, it then follows that $f(t, X_t)$ is itself a Dirichlet process. We also give a brief discussion later in this section of the possible applications of these results to one-dimensional diffusions.

Throughout this paper we assume the existence of a complete filtered probability space $(\Omega, \mathcal{F}, (\mathcal{F}_t)_{t \in \mathbb{R}_+}, \mathbb{P})$. The definition of quadratic variation used follows that of [17]. First, a (stochastic) partition $P$ of $\mathbb{R}_+$ is a sequence of stopping times $0 = \tau_0^P \leq \tau_1^P \leq \cdots \uparrow \infty$. Then for càdlàg processes $X, Y$ the approximation $[X, Y]^P$ to the quadratic covariation along a partition $P$ is

$$(2) \qquad [X, Y]_t^P \equiv \sum_{k=1}^{\infty} (X_{\tau_k^P \wedge t} - X_{\tau_{k-1}^P \wedge t})(Y_{\tau_k^P \wedge t} - Y_{\tau_{k-1}^P \wedge t}).$$

The quadratic covariation $[X, Y]$, if it exists, is defined to be the limit of $[X, Y]^P$ as the mesh $|P| \equiv \sup_{k \in \mathbb{N}} \|\tau_k^P - \tau_{k-1}^P\|_\infty$ goes to zero, with the topology of uniform convergence on compacts in probability (ucp for short).

$$[X, Y] = \lim_{|P| \to 0} [X, Y]^P \qquad \text{(ucp)}.$$

We also write $[X] \equiv [X, X]$ for the quadratic variation. If the quadratic variations and covariation of processes $X, Y$ all exist, then by the polarization identity $[X, Y] = ([X + Y] - [X] - [Y])/2$, the quadratic covariation is a difference of increasing processes and so has locally finite variation. By ucp convergence, the jumps of the quadratic covariation are $\Delta[X, Y] = \Delta X \Delta Y$, so its continuous part can be written as

$$[X, Y]_t^c \equiv [X, Y]_t - \sum_{s \leq t} \Delta X_s \Delta Y_s$$

and $[X]^c \equiv [X, X]^c$. A càdlàg process $X$ will be said to have *zero continuous quadratic variation* if its quadratic variation exists, and $[X]^c = 0$. Alternatively, for short, $X$ will be referred to as a z.c.q.v. process. Then the following definition of Dirichlet processes will be used.

DEFINITION 1.1. We say that a real valued process $X$ is a *Dirichlet process* if it has a decomposition $X = Y + V$ where $Y$ is a semimartingale and $V$ is a càdlàg adapted z.c.q.v. process.

We now state the first result which says that a locally Lipschitz continuous function of a semimartingale is a Dirichlet process. Such functions are differentiable almost everywhere, so we set $f'(x) \equiv \limsup_{h \to 0}(f(x+h) - f(x))/h$ which will be locally bounded and equal to the derivative of $f$ wherever it is differentiable.

THEOREM 1.2. *Let $X$ be a semimartingale and $f : \mathbb{R} \to \mathbb{R}$ be locally Lipschitz continuous. Then*

$$f(X_t) = \int_0^t f'(X_{s-}) \, dX_s + V_t,$$

*where $V$ has zero continuous quadratic variation.*

The proof of this is given in Section 3 and follows as a special case of the decomposition of functions of Dirichlet processes (Theorem 2.1).

For time-dependent functions, as well as requiring $f(t, x)$ to be locally Lipschitz continuous in $x$ with left and right derivatives everywhere, it will also be required that, locally, its variation in $t$ is integrable with respect to $x$. This leads us to look at the following classes of functions.

DEFINITION 1.3. We shall denote by $\mathcal{D}_0$ the set of functions $f : \mathbb{R}_+ \times \mathbb{R} \to \mathbb{R}$ such that:

- $f(t, x)$ is locally Lipschitz continuous in $x$ and càdlàg in $t$,

- for every $K_0 < K_1 \in \mathbb{R}$ and $T \in \mathbb{R}_+$ then

$$\int_{K_0}^{K_1} \int_0^T |d_t f(t,x)|\, dx < \infty.$$

If, furthermore, the left and right derivatives of $f(t,x)$ with respect to $x$ exist everywhere, then we write $f \in \mathcal{D}$.

As with the time-independent case above, the derivative of $f(t,x)$ with respect to $x$ need not exist everywhere, and the notation $D_x f(t,x)$ will be used to denote $\limsup_{h \to 0}(f(t,x+h) - f(t,x))/h$. Again, this will be locally bounded for any $f \in \mathcal{D}_0$ and equal to the partial derivative with respect to $x$ wherever it exists.

THEOREM 1.4.  *Let $X$ be a semimartingale and $f \in \mathcal{D}$. Then*

$$f(t, X_t) = \int_0^t D_x f(s, X_{s-})\, dX_s + V_t,$$

*where $V$ has zero continuous quadratic variation.*

In particular, this shows that $f(t, X_t)$ is a Dirichlet process. In Section 2 we state, and prove, a more general decomposition result which generalizes Theorem 1.4 to arbitrary Dirichlet processes. However, this result will also require $f(t,x)$ to be differentiable with respect to $x$ in an almost everywhere sense. Then, in Section 3 we show that if $X$ is a semimartingale then any function $f \in \mathcal{D}$ is indeed differentiable in the necessary "almost everywhere" sense, from which Theorem 1.4 follows.

Let us first discuss some possible applications of Theorems 1.2 and 1.4 to the theory of one-dimensional diffusions. Diffusions with a distributional drift have been studied in [8] and [9] via a generator $L$, written formally as

(3) $$Lf = \tfrac{1}{2}\sigma^2 f'' + b' f'.$$

Here $\sigma = \sigma(x)$ and $b = b(x)$ are continuous functions. The diffusion $X$ is then defined such that

$$f(X_t) - \int_0^t Lf(X_s)\, ds$$

is a local martingale for all functions $f$ in the domain of the generator $L$. If $b$ is not differentiable then (3) is understood only as a formal expression, and the full definition of $L$ and its domain are given in [8] and [9]. Let us consider the case where $b = \alpha \sigma^2/2$ for some $\alpha \in (0,1]$. Then $f$ is in the domain of the generator $L$, if $\sigma^{2\alpha} f'$ is continuously differentiable and

$$Lf = \tfrac{1}{2}\sigma^{2-2\alpha}(\sigma^{2\alpha} f')'.$$

In particular, $Lh = 0$ is solved by

$$h(x) = \int_0^x \sigma^{-2\alpha}(y)\,dy,$$

so $Y = h(X)$ is a local martingale. Then, $h^{-1}$ is continuously differentiable and it follows that $X = h^{-1}(Y)$ is a Dirichlet process.

We could consider extending this analysis to the case where $\sigma$ is merely bounded and measurable, such that $\sigma^{-2\alpha}$ is locally integrable. In that case $h^{-1}$ might not be differentiable, although it will be locally Lipschitz continuous. Then Theorem 1.2 shows that $X = h^{-1}(Y)$ will still be a Dirichlet process. Using Theorem 1.4, these ideas could be generalized to the case where $\sigma = \sigma(t, x)$ is time-dependent.

Another application of these results, which will be investigated in a future paper, is in obtaining generalizations of the backward Kolmogorov equation. Suppose, for the moment, that $X$ is a diffusion satisfying a stochastic differential equation of the form

$$(4) \qquad dX_t = \sigma(t, X_t)\,dW_t + b(t, X_t)\,dt$$

for a Brownian motion $W$. Given a twice continuously differentiable function $f(t, x)$, the backward equation says that $f(t, X_t)$ is a local martingale if

$$(5) \qquad \frac{\partial f}{\partial t} + \frac{1}{2}\sigma^2 \frac{\partial^2 f}{\partial x^2} + b\frac{\partial f}{\partial x} = 0,$$

which is a straightforward consequence of Itô's lemma. In particular, if $f$ is bounded and satisfies the boundary condition $f(T, x) = g(x)$, then (5) provides a sufficient condition for

$$(6) \qquad f(t, X_t) = \mathbb{E}[g(X_T)|\mathcal{F}_t]$$

to be satisfied for all $t < T$. Under sufficiently strong conditions for the coefficients $\sigma$ and $b$—such as Hölder continuity (see [11])—this can be used to prove uniqueness of solutions to (4). Now, suppose that $\sigma$, $b$ are not smooth (and more generally, could be distributions). Then requiring $f$ to be twice differentiable is too restrictive for the backward equation to be useful, and (4) can fail to have unique solutions. However, in many cases, it is sufficient to restrict to functions $f \in \mathcal{D}$. For example, if $X$ is a continuous and strong Markov martingale, then the results of [14] and [16] show that if $g$ is convex, then $f(t, x)$ satisfying (6) turns out to be convex in $x$ and decreasing in $t$.

As any local martingale with zero quadratic variation must be constant, Theorem 1.4 shows that $f(t, X_t)$ will be a local martingale if and only if

$$V_t = V_0 - \int_0^t D_x f(s, X_s) b(s, X_s)\,ds.$$

Using this idea, it is possible to derive generalizations of the backward equation which apply to nondifferentiable functions. We shall apply such methods in a future paper to obtain uniqueness results for time-inhomogeneous one-dimensional diffusions.

We end this section with a few remarks on Dirichlet and z.c.q.v. processes. First, the quadratic covariation $[X, Y]$ is easy to describe whenever either of $X$ or $Y$ has zero continuous quadratic variation.

LEMMA 1.5. *Let $X$ and $Y$ be càdlàg processes such that $X$ has zero continuous quadratic variation and $[Y]$ exists. Then the covariation $[X, Y]$ exists and satisfies $[X, Y]^c = 0$.*

The proof of this is given in Section 2. If $[X] = 0$ then this result reduces to the statement $[X, Y] = 0$, which is a simple consequence of the Cauchy–Schwarz inequality. One implication of Lemma 1.5 is that the sum of any two z.c.q.v. processes is itself a z.c.q.v. process, and it follows that the space of Dirichlet processes is closed under taking linear combinations.

Note that although the decomposition into a semimartingale and zero continuous quadratic variation process will not be unique, any Dirichlet process $X$ has the canonical decomposition

$$X = M + V, \tag{7}$$

where $M$ is a continuous local martingale and $V$ is a z.c.q.v. process with $V_0 = 0$. The existence of the decomposition follows from the existence for the case where $X$ is a semimartingale ([12], page 209 or [15], page 527). Uniqueness follows from the fact that any local martingale with zero quadratic variation is constant.

Alternatively, the following Doob–Meyer-style decomposition can be used and is a generalization of the canonical decomposition for special semimartingales.

LEMMA 1.6. *Let $X$ be a Dirichlet process such that $X_t^* \equiv \sup_{s \leq t} |X_s|$ is locally integrable. Then there exists a unique decomposition $X = M + V$ where $M$ is a local martingale and $V$ is a previsible z.c.q.v. process with $V_0 = 0$.*

PROOF. First, as every previsible local martingale is continuous, it follows that every previsible z.c.q.v. local martingale has zero quadratic variation and, therefore, is constant. So, the decomposition is unique.

Existence of the decomposition is trivial for local martingales, so, by decomposition (7), it is enough to consider the case where $X$ has zero continuous quadratic variation. Write $^p\Delta X$ for the previsible projection of the process $\Delta X$. Then Theorem 7.42 of [12] shows that there exists a local martingale $M$ such that $\Delta M = \Delta X - {}^p\Delta X$. By applying decomposition (7) to $M$, without loss of generality we may suppose that $M$ has zero continuous quadratic variation. Writing $V = X - M$ we see that $\Delta V = {}^p\Delta X$ is previsible, so $V$ is a previsible z.c.q.v. process. $\square$

**2. Functions of Dirichlet processes.** In this section we shall state and prove the most general decomposition result of this paper for functions of Dirichlet processes. As $f(t,x)$ will be required to be differentiable with respect to $x$ in an "almost everywhere" sense, we start by defining

$$\text{(8)} \qquad \text{diff}(f) = \{(t,x) \in \mathbb{R}_+ \times \mathbb{R} : f(t,x) \text{ is differentiable in } x\}.$$

We also define the subset of $\mathbb{R}_+ \times \mathbb{R}$ at which $f(t,x)$ is differentiable with respect to $x$ in a rather strong sense.

$$\text{(9)} \quad \text{diff}\,C(f) = \left\{(t,x) : \lim_{\substack{s \to t \\ y,z \to x}} (f(s,z) - f(s,y))/(z-y) \text{ exists}\right\} \subseteq \text{diff}(f).$$

Here, the limit is taken over all $s \in \mathbb{R}_+$ and $y, z \in \mathbb{R}$ with $y \neq z$. Alternatively, $\text{diff}\,C(f)$ is the set of points at which $D_x f$ is continuous.

We now state the decomposition result.

THEOREM 2.1. *Let $X = Y + Z$ where $Y$ is a semimartingale and $Z$ is a càdlàg adapted z.c.q.v. process. Let $f \in \mathcal{D}_0$ satisfy*

$$\text{(10)} \qquad \int 1_{\{(t,X_t) \notin \text{diff}(f)\}}\, d[X]^c_t = 0,$$

$$\text{(11)} \qquad \iint 1_{\{(t,x) \notin \text{diff}\,C(f), \mathbb{P}(X_t=x)>0\}} |d_t f(t,x)|\, dx = 0.$$

*Then*

$$\text{(12)} \qquad f(t, X_t) = \int_0^t D_x f(s, X_{s-})\, dY_s + V_t,$$

*where $V$ is a z.c.q.v. process.*

Equation (11) is trivially satisfied whenever $f$ is time independent, and, as will be shown in Lemma 3.2, it is always satisfied in the case where $X$ is a semimartingale.

The proof of Theorem 2.1 is given in this section. We start with a necessary and sufficient condition for a process to have zero continuous quadratic variation (Lemma 2.3). This result is used firstly to give a short proof of Lemma 1.5, and then applied to Theorem 2.1, the proof of which is split up into several lemmas.

Let us introduce some notation in order to simplify the formulas used in this section. For any process $X$ and stochastic partition $P$ of $\mathbb{R}_+$, we use $\delta_k^P X \equiv X_{\tau_k^P} - X_{\tau_{k-1}^P}$, so expression (2) can be written as

$$[X,Y]_t^P = \sum_{k>0} \delta_k^P X^t \delta_k^P Y^t.$$

Here, $X^t$ denotes the stopped process $X_s^t \equiv X_{s \wedge t}$. Now suppose that $X, Y$ are any càdlàg processes and $S \subseteq \mathbb{R}_+ \times \Omega$ is a jointly measurable set containing only finitely many times in each bounded time interval (restricting to any $\omega \in \Omega$). We shall make use of the following limit, in order to subtract out the discontinuities of $X$ and $Y$,

$$(13) \quad \lim_{|P| \to 0} \sum_{k=1}^{\infty} 1_{\{\rrbracket\tau_{k-1}^P, \tau_k^P\rrbracket \cap S \neq \varnothing\}} \delta_k^P X^t \delta_k^P Y^t = \sum_{s \leq t} 1_{\{s \in S\}} \Delta X_s \Delta Y_s.$$

This follows from the fact the left-hand side reduces to a finite sum with one term for each time in $\rrbracket 0, t \rrbracket \cap S$, and convergence is almost-surely uniform over finite time intervals. So, define $\mathcal{S}$ to be the collection of jointly measurable subsets of $\mathbb{R}_+ \times \Omega$ which contain only finitely many times in each bounded time interval (for each $\omega \in \Omega$). By the debut theorem ([5], IV.50 or [12], IV.1), this is the same as the sets which can be expressed as the union of graphs of a sequence of random times increasing to infinity,

$$(14) \quad \mathcal{S} = \left\{ \bigcup_{n=1}^{\infty} \llbracket \tau_n \rrbracket : \tau_n : \Omega \to \mathbb{R}_+ \cup \{\infty\} \text{ are measurable and } \tau_n \uparrow \infty \right\}.$$

For any partition $P$, $S \in \mathcal{S}$ and $t > 0$, we write $[P, S, t]$ to denote the (random) set of $k \in \mathbb{N}$ such that $\tau_k^P < t$ and $\rrbracket \tau_{k-1}^P, \tau_k^P \rrbracket \cap S$ is empty. Using this notation, we now give a sufficient condition for $[X, Y]^c = 0$ to be satisfied.

LEMMA 2.2. *Let $X$ and $Y$ be càdlàg adapted processes such that*

$$(15) \quad \inf_{S \in \mathcal{S}} \limsup_{|P| \to 0} \mathbb{P}\left( \sum_{k \in [P,S,t]} |\delta_k^P X \delta_k^P Y| > \varepsilon \right) = 0$$

*for all $t, \varepsilon > 0$. The limit is taken as $P$ ranges over the partitions of $\mathbb{R}_+$. Then the quadratic covariation $[X, Y]$ exists and $[X, Y]^c = 0$.*

PROOF. First, we note that for every $S \in \mathcal{S}$ and $t > 0$,

$$\sum_{s<t} |\Delta X_s \Delta Y_s| \leq \sum_{s \in S, s<t} |\Delta X_s \Delta Y_s| + \liminf_{|P| \to 0} \sum_{k \in [P,S,t]} |\delta_k^P X \delta_k^P Y|.$$

By (15), the right-hand side of this expression must, with probability 1, be finite for some $S \in \mathcal{S}$. Therefore, the locally-finite variation process $A_t =$

$\sum_{s\leq t}\Delta X_s\Delta Y_s$ is well defined. We show that $[X,Y]=A$. Consider the following identity:

$$[X,Y]_s^P - A_s = \sum_{k=1}^{\infty} 1_{\{\rrbracket\tau_{k-1}^P,\tau_k^P\rrbracket \cap S \neq \varnothing\}} \delta_k^P X^s \delta_k^P Y^s - \sum_{u \in S} \Delta X_u^s \Delta Y_u^s$$

$$+ \sum_{k=1}^{\infty} 1_{\{\rrbracket\tau_{k-1}^P,\tau_k^P\rrbracket \cap S = \varnothing\}} \delta_k^P X^s \delta_k^P Y^s - \sum_{u \notin S} \Delta A_u^s.$$

Limit (13) says that the first two terms on the right-hand side vanish as $|P|$ goes to zero (uniformly over all $s < t$), giving

$$\limsup_{|P|\to 0} \mathbb{P}\left(\sup_{s<t}|[X,Y]_s^P - A_s| \geq \varepsilon\right)$$

$$\leq \limsup_{|P|\to 0} \mathbb{P}\left(\sum_{k\in[P,S,t]} |\delta_k^P X \delta_k^P Y| \geq \varepsilon/2\right) + \mathbb{P}\left(\sum_{s\notin S, s<t} |\Delta A_s| \geq \varepsilon/2\right)$$

for all $t, \varepsilon > 0$. As $A$ is càdlàg and measurable, $S$ can be increased to include all the jump times of $A$ in the limit, so the last term on the right-hand side can be made arbitrarily small. Also, by the condition of the lemma, the first term can also be made as small as we like. $\square$

This leads to the following necessary and sufficient condition for a process to have zero continuous quadratic variation.

LEMMA 2.3. *Let $X$ be a càdlàg process. Then it has zero continuous quadratic variation if and only if*

(16) $$\inf_{S\in\mathcal{S}} \limsup_{|P|\to 0} \mathbb{P}\left(\sum_{k\in[P,S,t]} (\delta_k^P X)^2 > \varepsilon\right) = 0$$

*for all $t, \varepsilon > 0$.*

PROOF. If (16) is satisfied, then Lemma 2.2 with $Y = X$ gives the result. Conversely, suppose that $X$ has zero continuous quadratic variation and consider the following identity,

$$\sum_{k\in[P,S,t]} (\delta_k^P X)^2 = [X]_\tau^P - \sum_{s\leq\tau}(\Delta X_s)^2 + \sum_{s\notin S, s\leq\tau}(\Delta X_s)^2$$

$$+ \sum_{s\in S, s\leq\tau}(\Delta X_s)^2 - \sum_{\tau_k^P<t} 1_{\{\rrbracket\tau_{k-1}^P,\tau_k^P\rrbracket \cap S \neq \varnothing\}}(\delta_k^P X)^2.$$

Here, $\tau$ is the maximum of the stopping times $\tau_k^P$ satisfying $\tau_k^P < t$. As $X$ has zero continuous quadratic variation, the first two terms on the right-hand side converge to zero in probability as $|P|$ tends to 0. Also, limit (13) shows that the last two terms vanish, giving

$$\limsup_{|P|\to 0} \mathbb{P}\bigg(\sum_{k\in[P,S,t]} (\delta_k^P X)^2 > \varepsilon\bigg) \le \mathbb{P}\bigg(\sum_{s\notin S, s<t} (\Delta X_s)^2 \ge \varepsilon\bigg).$$

The result follows by noting that we can increase $S$ to include all the jump times of $X$ in the limit. $\square$

Lemma 1.5 follows as a simple consequence of Lemmas 2.2 and 2.3.

PROOF OF LEMMA 1.5. For $S \in \mathcal{S}$, $t > 0$ and partition $P$, the Cauchy–Schwarz inequality gives

$$\sum_{k\in[P,S,t]} |\delta_k^P X \delta_k^P Y| \le \bigg(\sum_{k\in[P,S,t]} (\delta_k^P X)^2\bigg)^{1/2} \bigg(\sum_{\tau_k^P < t} (\delta_k^P Y)^2\bigg)^{1/2}.$$

As the quadratic variation $[Y]$ is well defined we can take limits as $|P| \to 0$,

$$\limsup_{|P|\to 0} \mathbb{P}\bigg(\sum_{k\in[P,S,t]} |\delta_k^P X \delta_k^P Y| > \varepsilon\bigg)$$

$$\le \limsup_{|P|\to 0} \mathbb{P}\bigg([Y]_t \sum_{k\in[P,S,t]} (\delta_k^P X)^2 \ge \varepsilon^2\bigg)$$

$$\le \limsup_{|P|\to 0} \mathbb{P}\bigg(K \sum_{k\in[P,S,t]} (\delta_k^P X)^2 \ge \varepsilon^2\bigg) + \mathbb{P}([Y]_t > K)$$

for all $\varepsilon, K > 0$. As $X$ has zero continuous quadratic variation, Lemma 2.3 says that the first term on the right-hand side of this inequality goes to 0 if we take the infimum over all $S \in \mathcal{S}$. Then, taking the limit as $K \to \infty$, we see that the second term on the right-hand side also vanishes. So, the result follows from Lemma 2.2. $\square$

The remainder of this section is dedicated to proving Theorem 2.1. Let $V$ be the process appearing on the right-hand side of (12),

$$V_t \equiv f(t, X_t) - \int_0^t D_x f(s, X_{s-})\, dY_s.$$

It needs to be shown that this is a z.c.q.v. process, and the approach used is to split $\delta_k^P V$ into separate parts,

$$\delta_k^P V = (f(\tau_k^P, X_{\tau_k^P}) - f(\sigma, X_{\tau_k^P}) + f(\sigma, X_{\tau_{k-1}^P}) - f(\tau_{k-1}^P, X_{\tau_{k-1}^P}))$$

$$\text{(17)} \qquad + \left( \zeta_\sigma \delta_k^P X - \int_{\tau_{k-1}^P}^{\tau_k^P} D_x f(t, X_{t-})\, dY_t \right)$$

$$+ (f(\sigma, X_{\tau_k^P}) - f(\sigma, X_{\tau_{k-1}^P}) - \zeta_\sigma \delta_k^P X).$$

Here, $\sigma$ is a suitably chosen stopping time in the interval $[\tau_{k-1}^P, \tau_k^P]$ and $\zeta$ is a simple previsible process which, by definition, are linear combinations of processes of the form $A1_{\{t > \tau\}}$ for stopping times $\tau$ and bounded $\mathcal{F}_\tau$-measurable random variables $A$.

Using Lemma 2.3, we show that the contribution of each of the three terms on the right-hand side of (17) to the continuous part of the quadratic variation of $V$ can be made arbitrarily small (by making suitable choices of $\sigma$ and $\zeta$).

We start by showing that the contribution to the continuous part of the quadratic variation coming from the first term on the right-hand side of (17) is zero. The idea is to smooth out the time increments of $f$ by making use of the following identity:

$$\text{(18)} \qquad g(y) = \frac{1}{a} \int_{y-a}^{y} ((a - y + x)g'(x) + g(x))\, dx,$$

which is an application of integration by parts and applies for every absolutely continuous function $g$ and every $a > 0$.

LEMMA 2.4. *Let $X$ be a càdlàg process and let $f \in \mathcal{D}_0$ satisfy (11). Then for any $t > 0$,*

$$\operatorname*{ess\,inf}_{S \in \mathcal{S}} \limsup_{|P| \to 0} \sum_{k \in [P, S, t]} \sup_{s \in [\tau_{k-1}^P, \tau_k^P]} (f(\tau_k^P, X_s) - f(\tau_{k-1}^P, X_s))^2 = 0.$$

PROOF. For any $u < v \in \mathbb{R}_+$ and $x \in \mathbb{R}$, we use the notation

$$\delta_{u,v} f(x) \equiv f(v, x) - f(u, x).$$

Then for any $a > 0$, substituting $g(x) = (\delta_{u,v} f(x))^2$ into (18) gives

$$(\delta_{u,v} f(y))^2 = \frac{1}{a} \int_{y-a}^{y} (2(a - y + x)(\delta_{u,v} D_x f(x))\delta_{u,v} f(x) + (\delta_{u,v} f(x))^2)\, dx$$

$$= \frac{1}{a} \int_{y-a}^{y} \int_{u}^{v} (2(a - y + x)\delta_{u,v} D_x f(x) + \delta_{u,v} f(x))\, d_t f(t, x)\, dx.$$

For any $S \in \mathcal{S}$, it follows that if $h_a^{P,S}(u, x)$ is the (random) function

$$h_a^{P,S}(u, x) = \frac{1}{a} \sum_{k=1}^{\infty} 1_{\{\tau_{k-1}^P < u \leq \tau_k^P\}} 1_{\{]\!]\tau_{k-1}^P, \tau_k^P]\!] \cap S = \varnothing\}} \sup_{s \in [\tau_{k-1}^P, \tau_k^P]} 1_{\{x \in (X_s - a, X_s)\}}$$

$$\times |2(a - X_s + x)\delta_{\tau_{k-1}^P, \tau_k^P} D_x f(x) + \delta_{\tau_{k-1}^P, \tau_k^P} f(x)|,$$

then
$$A_S^P \equiv \sum_{k\in[P,S,t]} \sup_{s\in[\tau_k^P,\tau_{k-1}^P]} (\delta_{\tau_{k-1}^P,\tau_k^P} f(X_s))^2 \leq \int_{-\infty}^{\infty}\int_0^t h_a^{P,S}(s,x)|d_s f(s,x)|\,dx.$$

Without loss of generality, we can assume that $f(t,x)$ is Lipschitz continuous in $x$ with coefficient $K$, in which case
$$\limsup_{|P|\to 0} |h_a^{P,S}(s,x)| \leq 1_{\{s\notin S\}} g_a(s,x),$$
$$g_a(s,x) \equiv 1_{\{X_{s-}\wedge X_s - a \leq x \leq X_{s-}\vee X_s\}}$$
$$\times (4K1_{\{(s,x)\notin \mathrm{diff}\,C(f)\}} + a^{-1}|\Delta_s f(s,x)|),$$
where $\Delta_s f(s,x) \equiv f(s,x) - f(s-,x)$. So, by bounded convergence,
$$\limsup_{|P|\to 0} A_S^P \leq \int_{-\infty}^{\infty}\int_0^t 1_{\{s\notin S\}} g_a(s,x)|d_s f(s,x)|\,dx.$$

As $S \in \mathcal{S}$ can be increased to include (in the limit) all the times at which either $f(s,x)$ or $X_s$ is not continuous,
$$\operatorname*{ess\,inf}_{S\in\mathcal{S}} \limsup_{|P|\to 0} A_S^P \leq 2K\int_{-\infty}^{\infty}\int_0^t 1_{\{X_s-a\leq x\leq X_s,(s,x)\notin \mathrm{diff}\,C(f)\}}|d_s f(s,x)|\,dx.$$

Also, $a$ can be chosen arbitrarily small,
$$\operatorname*{ess\,inf}_{S\in\mathcal{S}} \limsup_{|P|\to 0} A_S^P \leq 2K\int_{-\infty}^{\infty}\int_0^t 1_{\{x=X_s,(s,x)\notin \mathrm{diff}\,C(f)\}}|d_s f(s,x)|\,dx.$$

Finally, (11) shows that the right-hand side has zero expectation, so it must almost surely be equal to 0. □

We now bound the contribution to the continuous part of the quadratic variation of $V$ coming from the second term on the right-hand side of (17). The previsible process $\eta_s$ below will be chosen to be equal to $D_x f(s, X_{s-})$.

LEMMA 2.5. *Let $X = Y + Z$ where $Y$ is a semimartingale and $Z$ is a càdlàg adapted process with zero continuous quadratic variation. Given any uniformly bounded and previsible process $\eta$ and simple previsible process $\zeta$ set,*
$$B_k^P \equiv \sup_{s\in[\tau_{k-1}^P,\tau_k^P]} \left|\zeta_s \delta_k^P X - \int_{\tau_{k-1}^P}^{\tau_k^P} \eta_u\,dY_u\right|$$
*for all partitions $P$ of $\mathbb{R}_+$. Then*
$$\inf_{S\in\mathcal{S}} \limsup_{|P|\to 0} \mathbb{P}\left(\sum_{k\in[P,S,t]} (B_k^P)^2 \geq \varepsilon\right) \leq \mathbb{P}\left(\int_0^t (\zeta - \eta)^2\,d[Y] \geq \varepsilon\right)$$

*for all $t, \varepsilon > 0$.*

PROOF. First, as $\zeta$ is simple previsible, it is piecewise constant and there are only finitely many times at which it is not continuous. So, we can restrict to those $S \in \mathcal{S}$ which contain all the discontinuity times of $\zeta$. In that case, for any $k \in [P, S, t]$ and $s \in (\tau_{k-1}^P, \tau_k^P]$, we have $\zeta_s = \zeta_{\tau_k^P}$. So, for $k \in [P, S, t]$,

$$\pm B_k^P = \zeta_{\tau_k^P} \delta_k^P Z + \left( \zeta_{\tau_k^P} \delta_k^P Y - \int_{\tau_{k-1}^P}^{\tau_k^P} \eta_s \, dY_s \right)$$
$$= \zeta_{\tau_k^P} \delta_k^P Z + \delta_k^P U,$$

where $U$ is the process $U = \int (\zeta - \eta) \, dY$. Then the triangle inequality gives

$$(19) \quad \left( \sum_{k \in [P,S,t]} (B_k^P)^2 \right)^{1/2} \leq K \left( \sum_{k \in [P,S,t]} (\delta_k^P Z)^2 \right)^{1/2} + \left( \sum_{k=1}^{\infty} (\delta_k^P U^t)^2 \right)^{1/2},$$

where $K$ is any upper bound for $|\zeta|$. As $Z$ has zero continuous quadratic variation, Lemma 2.3 gives

$$\inf_{S \in \mathcal{S}} \limsup_{|P| \to 0} \mathbb{P} \left( \sum_{k \in [P,S,t]} (\delta_k^P Z)^2 \geq \varepsilon \right) = 0$$

for every $\varepsilon > 0$. Finally, using the definition of quadratic variation, the last term on the right-hand side of inequality (19) will converge in probability to $[U]_t = \int_0^t (\zeta - \eta)^2 \, d[Y]$ as $|P| \to 0$, giving the result. $\square$

We now turn to the third term on the right-hand side of (17). This will require making a suitable choice for $\sigma \in [\tau_{k-1}^P, \tau_k^P]$. More precisely, for every partition $P$, we will choose stopping times $(\sigma_k^P)_{k \in \mathbb{N}}$ satisfying

$$(20) \quad \begin{aligned} & \tau_{k-1}^P \leq \sigma_k^P \leq \tau_k^P, \\ & \sigma_k^P > \tau_{k-1}^P \quad \text{whenever } \tau_k^P > \tau_{k-1}^P, \end{aligned}$$

for each $k$. Once these times have been chosen, they define a new partition $\tilde{P}$ given by

$$(21) \quad \tau_k^{\tilde{P}} = \begin{cases} \tau_{k/2}^P, & \text{if } k \text{ is even}, \\ \sigma_{(k+1)/2}^P, & \text{if } k \text{ is odd}. \end{cases}$$

The choice of $\sigma_k^P$ will be made with the help of the following lemma, the proof of which makes use of the optional section theorem ([5], IV.84 or [12], Theorem 4.7).

LEMMA 2.6. *Let $X$ be a Dirichlet process, and $\xi$ be any nonnegative optional process uniformly bounded by some $K \in \mathbb{R}_+$. For every partition $P$ set*

$$D_k^P \equiv \xi_{\sigma_k^P}((X_{\tau_k^P} - X_{\sigma_k^P})^2 + (X_{\sigma_k^P} - X_{\tau_{k-1}^P})^2).$$

*Then for every $\delta > 0$ we can choose the stopping times $\sigma_k^P$ satisfying inequalities (20) such that*

$$\inf_{S \in \mathcal{S}} \limsup_{|P| \to 0} \mathbb{P}\bigg(\sum_{k \in [P,S,t]} D_k^P \geq \varepsilon\bigg) \leq \mathbb{P}\bigg(\int_0^t (K \mathbf{1}_{\{|\xi_s| > \delta\}} + \delta) \, d[X]_s^c \geq \varepsilon\bigg)$$

*for all $t, \varepsilon > 0$.*

PROOF. First, by decomposition (7), write $X = Y + Z$ for a continuous local martingale $Y$ and z.c.q.v. process $Z$. Let us set

$$A_k^P \equiv \xi_{\sigma_k^P}((Y_{\tau_k^P} - Y_{\sigma_k^P})^2 + (Y_{\sigma_k^P} - Y_{\tau_{k-1}^P})^2),$$
$$B_k^P \equiv \xi_{\sigma_k^P}((Z_{\tau_k^P} - Z_{\sigma_k^P})^2 + (Z_{\sigma_k^P} - Z_{\tau_{k-1}^P})^2).$$

Choosing any $S \in \mathcal{S}$ the triangle inequality gives

$$\bigg(\sum_{k \in [P,S,t]} D_k^P\bigg)^{1/2} \leq \bigg(\sum_{\tau_k^P < t} A_k^P\bigg)^{1/2} + \bigg(\sum_{k \in [P,S,t]} B_k^P\bigg)^{1/2}$$
$$\leq \bigg(\sum_{\tau_k^P < t} A_k^P\bigg)^{1/2} + \bigg(K \sum_{k \in [\tilde{P},S,t]} (\delta_k^{\tilde{P}} Z)^2\bigg)^{1/2},$$

where $\tilde{P}$ is the partition defined by (21). If we choose any $\varepsilon' < \varepsilon$ and set $\varepsilon'' = K^{-1}(\sqrt{\varepsilon} - \sqrt{\varepsilon'})^2$, this gives

$$\mathbb{P}\bigg(\sum_{k \in [P,S,t]} D_k^P \geq \varepsilon\bigg) \leq \mathbb{P}\bigg(\sum_{\tau_k^P < t} A_k^P \geq \varepsilon'\bigg) + \mathbb{P}\bigg(\sum_{k \in [\tilde{P},S,t]} (\delta_k^{\tilde{P}} Z)^2 \geq \varepsilon''\bigg).$$

As $Z$ has zero continuous quadratic variation, Lemma 2.3 says that the second term on the right-hand side vanishes if we let $|P|$ go to zero and take the infimum over all $S \in \mathcal{S}$,

(22) $\quad \inf_{S \in \mathcal{S}} \limsup_{|P| \to 0} \mathbb{P}\bigg(\sum_{k \in [P,S,t]} D_k^P \geq \varepsilon\bigg) \leq \limsup_{|P| \to 0} \mathbb{P}\bigg(\sum_{\tau_k^P < t} A_k^P \geq \varepsilon'\bigg).$

This simplifies the problem to the case of a continuous local martingale.

We now make a choice for the stopping times $\sigma_k^P$. For any partition $P$ and $k \in \mathbb{N}$, the set of times $s \in (\tau_{k-1}^P, \tau_k^P]$ such that $\xi_s \leq \delta$ is optional. So, by the optional section theorem the stopping time $\sigma_k^P$ can be chosen such that inequalities (20) are satisfied, $\xi_{\sigma_k^P} \leq \delta$ whenever $\sigma_k^P < \tau_k^P$ and

$$\mathbb{P}(\sigma_k^P < \tau_k^P) \geq \mathbb{P}(\exists s \in (\tau_{k-1}^P, \tau_k^P) \text{ s.t. } \xi_s \leq \delta) - 2^{-k}|P|.$$

It follows that

$$\mathbb{P}(\xi_{\sigma_k^P} > \delta) \leq \mathbb{P}(\forall s \in (\tau_{k-1}^P, \tau_k^P], \xi_s > \delta) + 2^{-k}|P|.$$

Also, by the debut theorem, we can define the stopping times

$$\tilde{\sigma}_k^P = \inf\{s \in (\tau_{k-1}^P, \tau_k^P] : \xi_s \leq \delta\} \cup \{\tau_k^P\}.$$

By the choice of $\sigma_k^P$ and $\tilde{\sigma}_k^P$, the following holds outside of a set of probability at most $2^{-k}|P|$:

$$A_k^P \leq \delta((Y_{\tau_k^P} - Y_{\sigma_k^P})^2 + (Y_{\sigma_k^P} - Y_{\tau_{k-1}^P})^2) + K 1_{\{\sigma_k^P = \tau_k^P\}}(Y_{\tau_k^P} - Y_{\tau_{k-1}^P})^2$$

$$\leq \delta((\delta_{2k}^{\tilde{P}} Y)^2 + (\delta_{2k-1}^{\tilde{P}} Y)^2) + K(Y_{\tilde{\sigma}_k^P} - Y_{\tau_{k-1}^P})^2$$

$$= \delta(\delta_{2k}^{\tilde{P}} Y)^2 + \delta(\delta_{2k-1}^{\tilde{P}} Y)^2 + 2K \int_{\tau_{k-1}^P}^{\tilde{\sigma}_k^P} (Y_s - Y_{\tau_{k-1}^P}) \, dY_s + K \int_{\tau_{k-1}^P}^{\tilde{\sigma}_k^P} d[Y]_s,$$

where $\tilde{P}$ is the partition defined by (21). Noting that $\xi_s > \delta$ whenever $s \in (\tau_{k-1}^P, \tilde{\sigma}_k^P)$, this inequality gives

$$(23) \quad A_k^P \leq \delta(\delta_{2k}^{\tilde{P}} Y)^2 + \delta(\delta_{2k-1}^{\tilde{P}} Y)^2 + K \int_{\tau_{k-1}^P}^{\tau_k^P} 1_{\{\xi_s > \delta\}} d[Y]_s + 2K \int_{\tau_{k-1}^P}^{\tau_k^P} \alpha_s^P \, dY_s$$

outside of a set with probability at most $2^{-k}|P|$ and with

$$\alpha_s^P \equiv \sum_{k=1}^{\infty} 1_{\{s \in (\tau_{k-1}^P, \tilde{\sigma}_k^P]\}}(Y_s - Y_{\tau_{k-1}^P}).$$

The continuity of $Y$ implies that $\alpha^P \to 0$ as $|P| \to 0$, so bounded convergence for stochastic integration gives

$$\sup_{s<t} \left| \int_0^s \alpha_u^P \, dY_u \right| \to 0$$

in probability as $|P| \to 0$. Summing inequality (23) over $k$ and taking the limit as $|P| \to 0$ gives

$$\limsup_{|P| \to 0} \mathbb{P}\left(\sum_{\tau_k^P < t} A_k^P \geq \varepsilon'\right)$$

$$\leq \limsup_{|P| \to 0} \mathbb{P}\left(\delta \sum_{\tau_k^{\tilde{P}} < t} (\delta_k^{\tilde{P}} Y)^2 + K \int_0^t 1_{\{\xi_s > \delta\}} d[Y]_s > \tilde{\varepsilon}\right) + \limsup_{|P| \to 0} \sum_{k=0}^{\infty} 2^{-k}|P|$$

$$\leq \mathbb{P}\left(\int_0^t (\delta + K 1_{\{\xi_s > \delta\}}) \, d[Y]_s \geq \tilde{\varepsilon}\right),$$

where $\tilde{\varepsilon}$ is any real number in the range $0 < \tilde{\varepsilon} < \varepsilon'$. The result now follows from combining this with inequality (22) and letting $\tilde{\varepsilon}$ increase to $\varepsilon$. □

We use Lemma 2.6 to bound the contribution to the continuous part of the quadratic variation of $V$ coming from the third term on the right-hand side of (17).

LEMMA 2.7. *Let $X$ be a Dirichlet process and $f:\mathbb{R}_+ \times \mathbb{R} \to \mathbb{R}$ be càdlàg in $t$ and Lipschitz continuous in $x$. Choosing any bounded optional process $\zeta$ and any $h > 0$, set*

$$\xi_s \equiv \sup_{0 < |a| \leq h} |(f(s, X_s + a) - f(s, X_s))/a - \zeta_s|.$$

*Also, for every partition $P$, set*

$$C_k^P \equiv f(\sigma_k^P, X_{\tau_k^P}) - f(\sigma_k^P, X_{\tau_{k-1}^P}) - \zeta_{\sigma_k^P} \delta_k^P X.$$

*Then for any $\delta > 0$ the stopping times $\sigma_k^P$ satisfying inequalities (20) can be chosen such that*

$$\inf_{S \in \mathcal{S}} \limsup_{|P| \to 0} \mathbb{P}\left(\sum_{k \in [P, S, t]} (C_k^P)^2 \geq \varepsilon\right) \leq \mathbb{P}\left(2 \int_0^t (1_{\{\xi_s > \delta\}} K^2 + \delta^2) \, d[X]_s^c \geq \varepsilon\right)$$

*for all $t, \varepsilon > 0$ where $K \in \mathbb{R}$ is any upper bound for $\xi$.*

PROOF. First note that we can restrict $a$ to the rational numbers in the definition of $\xi$, so it is the supremum of a countable set of optional processes and therefore is itself optional.

For every partition $P$, set

$$a_k^P \equiv X_{\tau_k^P} - X_{\sigma_k^P}, \qquad b_k^P \equiv X_{\tau_{k-1}^P} - X_{\sigma_k^P}.$$

Then we can rewrite $C_k^P$ as

$$C_k^P = 1_{\{a_k^P \neq 0\}}((f(\sigma_k^P, X_{\sigma_k^P} + a_k^P) - f(\sigma_k^P, X_{\sigma_k^P}))/a_k^P - \zeta_{\sigma_k^P})a_k^P$$
$$- 1_{\{b_k^P \neq 0\}}((f(\sigma_k^P, X_{\sigma_k^P} + b_k^P) - f(\sigma_k^P, X_{\sigma_k^P}))/b_k^P - \zeta_{\sigma_k^P})b_k^P.$$

In particular, if $|a_k^P|$ and $|b_k^P|$ are both smaller than $h$, then

$$|C_k^P| \leq \xi_{\sigma_k^P}(|X_{\tau_k^P} - X_{\sigma_k^P}| + |X_{\sigma_k^P} - X_{\tau_{k-1}^P}|)$$

and so

(24) $\qquad (C_k^P)^2 \leq B_k^P \equiv 2\xi_{\sigma_k^P}^2((X_{\tau_k^P} - X_{\sigma_k^P})^2 + (X_{\sigma_k^P} - X_{\tau_{k-1}^P})^2).$

So, if we let $S \in \mathcal{S}$ include all the times $s$ for which $|\Delta X_s| \geq h$, then inequality (24) will hold whenever $]\!]\tau_{k-1}^P, \tau_k^P]\!] \cap S = \varnothing$ and $\tau_k^P < t$ for all fine enough partitions $P$. Therefore,

$$\mathbb{P}\bigg(\sum_{k \in [P,S,t]} (C_k^P)^2 \geq \varepsilon\bigg) \leq \mathbb{P}\bigg(\sum_{k \in [P,S,t]} B_k^P \geq \varepsilon\bigg)$$

in the limit as $|P| \to 0$. The result now follows by applying Lemma 2.6 with $2\xi^2$ in place of $\xi$, $2K^2$ in place of $K$, and $2\delta^2$ in place of $\delta$. $\square$

Finally, for this section, we put together the results of Lemmas 2.4, 2.5 and 2.7 to prove Theorem 2.1.

PROOF OF THEOREM 2.1. By the condition of the theorem, $X = Y + Z$ for semimartingale $Y$ and z.c.q.v. process $Z$. Using decomposition (7) we may suppose that $Y$ is continuous, so $[Y] = [X]^c$. It needs to be shown that $V$ defined by (12) has zero continuous quadratic variation. By localization, we may assume that $f(t, x)$ is Lipschitz continuous in $x$ with coefficient $L$, rather than just locally Lipschitz.

Let $\eta$ be the previsible process $\eta_s = D_x f(s, X_{s-})$, which is uniformly bounded by $L$. Also pick any simple previsible process $\zeta$ such that $|\zeta| \leq L$. For any $h > 0$, set

$$\xi_s^h \equiv \sup_{0 < |a| \leq h} |(f(s, X_s + a) - f(s, X_s))/a - \zeta_s|,$$

which is bounded by $2L$. Supposing that for every partition $P$ stopping times $\sigma_k^P$ satisfying inequalities (20) have been chosen, (17) allows us to write

$$\delta_k^P V = A_k^P + B_k^P + C_k^P$$

with

$$A_k^P = f(\tau_k^P, X_{\tau_k^P}) - f(\sigma_k^P, X_{\tau_k^P}) + f(\sigma_k^P, X_{\tau_{k-1}^P}) - f(\tau_{k-1}^P, X_{\tau_{k-1}^P}),$$

$$B_k^P = \zeta_{\sigma_k^P} \delta_k^P X - \int_{\tau_{k-1}^P}^{\tau_k^P} \eta_s \, dY_s,$$

$$C_k^P = f(\sigma_k^P, X_{\tau_k^P}) - f(\sigma_k^P, X_{\tau_{k-1}^P}) - \zeta_{\sigma_k^P} \delta_k^P X,$$

where $\sigma_k^P$ are stopping times satisfying inequalities (20). In particular,

(25) $$(\delta_k^P V)^2 \leq 3(A_k^P)^2 + 3(B_k^P)^2 + 3(C_k^P)^2.$$

If $\tilde{P}$ is the partition defined by (21), then Lemma 2.4 with $\tilde{P}$ in place of $P$ gives

$$\operatorname*{ess\,inf}_{S \in \mathcal{S}} \limsup_{|P| \to 0} \sum_{k \in [P,S,t]} (A_k^P)^2 = 0$$

for all $t > 0$. So, by applying Lemmas 2.5 and 2.7, respectively, to the second and third terms on the right-hand side of (25), for any $\delta > 0$, the stopping times $\sigma_k^P$ can be chosen so that

$$\inf_{S \in \mathcal{S}} \limsup_{|P| \to 0} \mathbb{P}\bigg(\sum_{k \in [P,S,t]} (\delta_k^P V)^2 \geq \varepsilon\bigg)$$

(26)
$$\leq \mathbb{P}\bigg(\int_0^t (\zeta - \eta)^2 \, d[X]_t^c \geq \varepsilon/3\bigg)$$

$$+ \mathbb{P}\bigg(2 \int_0^t (1_{\{\xi_s^h > \delta\}} 4L^2 + \delta^2) \, d[X]_s^c \geq \varepsilon/3\bigg)$$

for any $\varepsilon > 0$. Also, whenever $(s, X_s) \in \mathrm{diff}(f)$ then the definition of $\xi^h$ gives

$$\xi_s^h \to |D_x f(s, X_s) - \zeta_s|$$

as $h \to 0$. By (10), this limit holds almost everywhere with respect to the measure $\int_0^t \cdot d[X]^c$. Combining this with the inequality $1_{\{\xi^h > \delta\}} \leq \delta^{-2}(\xi^h)^2$, we can take limits as $h \to 0$ in inequality (26),

$$\inf_{S \in \mathcal{S}} \limsup_{|P| \to 0} \mathbb{P}\bigg(\sum_{k \in [P,S,t]} (\delta_k^P V)^2 \geq \varepsilon\bigg)$$

(27)
$$\leq \mathbb{P}\bigg(\int_0^t (\zeta - \eta)^2 \, d[X]^c \geq \varepsilon/3\bigg)$$

$$+ \mathbb{P}\bigg(2 \int_0^t (4\delta^{-2} L^2 (D_x f(s, X_s) - \zeta_s)^2 + \delta^2) \, d[X]_s^c \geq \varepsilon/3\bigg).$$

As the simple previsible processes generate the previsible $\sigma$-algebra, the monotone class lemma shows that there exists a sequence of simple previsible processes $\zeta^n$ satisfying

$$\mathbb{P}\bigg(\int_0^t (\zeta_s^n - \eta_s)^2 \, d[X]_s^c \geq \varepsilon\bigg) \to 0$$

as $n \to \infty$ for every $\varepsilon > 0$. Furthermore, if $\eta$ is bounded by $L$, then $\zeta^n$ can also be chosen to be bounded by $L$. So, we can substitute $\zeta^n$ for $\zeta$ in the right-hand side of inequality (27) and take limits

$$\inf_{S \in \mathcal{S}} \limsup_{|P| \to 0} \mathbb{P}\bigg(\sum_{k \in [P,S,t]} (\delta_k^P V)^2 \geq \varepsilon\bigg)$$

$$\leq \mathbb{P}\bigg(2 \int_0^t (4\delta^{-2} L^2 (D_x f(s, X_s) - \eta_s)^2 + \delta^2) \, d[X]_t^c \geq \varepsilon/3\bigg)$$

$$= \mathbb{P}(2\delta^2 [X]_t^c \geq \varepsilon/3).$$

This last equality holds because $\eta_s = D_x f(s, X_s)$ whenever $\Delta X_s = 0$. The result now follows by letting $\delta$ decrease to 0 and applying Lemma 2.3. □

**3. Functions of semimartingales.** In this section, the decomposition result Theorem 2.1 is applied to the case where $X$ is a semimartingale. Using Lemma A.3 for the "almost everywhere" differentiability of functions in $\mathcal{D}$, it is shown that (11) is automatically satisfied, and (10) is satisfied for every $f \in \mathcal{D}$. Theorems 1.2 and 1.4 then follow.

We start with the following simple result, which allows us to represent the marginal distributions of a semimartingale by a function $C \in \mathcal{D}$.

LEMMA 3.1. *Let $X$ be an càdlàg adapted process which decomposes as $X = M + A$ for a martingale $M$ and integrable process $A$ with integrable variation over each finite time interval. Define the function $C : \mathbb{R}_+ \times \mathbb{R} \to \mathbb{R}$ by $C(t, x) \equiv \mathbb{E}[(X_t - x)_+]$. Then $C(t, x)$ is convex in $x$, càdlàg in $t$ and for every $x \in \mathbb{R}$,*

$$C(t, x) + \mathbb{E}\left[\int_0^t |dA_s|\right]$$

*is increasing in $t$. In particular, $C \in \mathcal{D}$.*

PROOF. First, $(X_t - x)_+$ is convex in $x$, so by the linearity of expectations, $C(t, x)$ will also be convex in $x$. Also, from the decomposition of $X$ we see that $\{X_t : t \leq T\}$ is uniformly integrable for every $T > 0$. Therefore, as $(X_t - x)_+$ is càdlàg in $t$ we see that $C(t, x)$ will also be càdlàg.

Let us now set $f(t) \equiv \mathbb{E}[\int_0^t |dA_s|]$. Then for every $s < t$, Jensen's inequality $\mathbb{E}[(M_t + A_s - x)_+] \geq \mathbb{E}[(M_s + A_s - x)_+]$ gives

$$C(t, x) = \mathbb{E}[(M_t + A_t - x)_+] \geq \mathbb{E}[(M_t + A_s - x)_+] - \mathbb{E}[(A_t - A_s)_-]$$
$$\geq \mathbb{E}[(M_s + A_s - x)_+] - f(t) + f(s) = C(s, x) - f(t) + f(s).$$

So $C(t, x) + f(t)$ is increasing in $t$.

It only remains to show that $C \in \mathcal{D}$. First, the convexity in $x$ shows that $C(t, x)$ is locally Lipschitz continuous with left and right derivatives in $x$. Secondly, $C(t, x)$ can be expressed as the sum of $C(t, x) + f(t)$ and $-f(t)$, which are monotonic in $t$. So its variation satisfies

$$\int_0^T |d_t C(t, x)| \leq C(T, x) + 2f(T),$$

which is locally bounded. □

Equality (11) follows easily for semimartingales.

LEMMA 3.2. *Let $X$ be a semimartingale and $f \in \mathcal{D}_0$. Then*

$$\iint 1_{\{(t,x)\notin \operatorname{diff} C(f), \mathbb{P}(X_t=x)>0\}} |d_t f(t,x)| \, dx = 0.$$

PROOF. As $X$ is a semimartingale it decomposes as $X = M + A$ for a local martingale $M$ and finite variation process $A$. By pre-localization, we only need to consider the case where $\sup_{t\geq 0} |X_t|$ is integrable and, therefore, $A$ has locally integrable variation. Then, by localization, we may suppose that $A$ has integrable variation, and $M$ is a uniformly integrable martingale.

We now set $C(t,x) \equiv \mathbb{E}[(X_t - x)_+]$. Then, letting $D_x^- C$, $D_x^+ C$ be its left and right derivatives in $x$, respectively,

$$\mathbb{P}(X_t = x) = D_x^+ C(t,x) - D_x^- C(t,x).$$

Therefore,

$$\iint 1_{\{(t,x)\notin \operatorname{diff} C(f), \mathbb{P}(X_t=x)>0\}} |d_t f(t,x)| \, dx$$
$$\leq \iint 1_{\{D_x^+ C(t,x) \neq D_x^- C(t,x)\}} |d_t f(t,x)| \, dx$$
$$= \iint 1_{\{(t,x)\notin \operatorname{diff}(C)\}} |d_t f(t,x)| \, dx.$$

However, Lemma 3.1 says that $C \in \mathcal{D}$, so by Lemma A.3, the right-hand side of the above equality is 0. □

In order to complete the proof of Theorems 1.2 and 1.4 it is necessary to show that equality (10) is satisfied. The following identity, which follows from Itô's lemma, will be be used to this end.

LEMMA 3.3. *Let $X$ be a càdlàg adapted process which decomposes as $X = M + A$ for a martingale $M$ and càdlàg integrable process $A$ with integrable variation over finite time intervals. Set $C(t,x) \equiv \mathbb{E}[(X_t - x)_+]$ so that, by Lemma 3.1, $C \in \mathcal{D}$.*

*Then, for any nonnegative and measurable $\theta : \mathbb{R}_+ \times \mathbb{R} \to \mathbb{R}$ with bounded support,*

$$\iint \theta \, d_t C \, dx = \frac{1}{2} \mathbb{E}\left[\int_0^\infty \theta(t, X_t) \, d[X]_t^c\right] + \mathbb{E}\left[\int_0^\infty \int_{-\infty}^{X_{t-}} \theta(t,y) \, dy \, dA_t\right]$$
$$\text{(28)} \qquad + \mathbb{E}\left[\sum_{t\in\mathbb{R}_+} \int_{X_{t-}}^{X_t} (X_t - x)\theta(t,x) \, dx\right].$$

PROOF. It is enough to consider the case where $\theta(t,x)$ is nonnegative, twice continuously differentiable in $x$ and once in $t$, and with compact support in $(0,\infty) \times \mathbb{R}$. The general case follows from the monotone class lemma. So suppose that $\theta$ satisfies these properties and define $f : \mathbb{R}_+ \times \mathbb{R} \to \mathbb{R}$ by

$$f(t,x) = \int \theta(t,y)(x-y)_+ \, dy,$$

which is twice continuously differentiable in $x$ with $D_{xx}f = \theta$. Also, as $\theta$ has compact support, $f$ has bounded derivatives and $0 \leq f(t,x) \leq K(1+|x|)$ for some constant $K$. Then Itô's lemma gives

$$f(t,X_t) = \int_0^t D_x f(s, X_{s-}) \, dX_s + \frac{1}{2} \int_0^t \theta(s, X_s) \, d[X]_s^c$$
$$+ \int_0^t D_t f(s, X_s) \, ds + \sum_{s \leq t} \int_{X_{s-}}^{X_s} (X_s - x)\theta(s,x) \, dx.$$

As $\int_0^t D_x f(s, X_{s-}) \, dM_s$ is a local martingale, there exist stopping times $T_n \uparrow \infty$ such that

$$\mathbb{E}\left[\int_0^{t \wedge T_n} D_x f(s, X_{s-}) \, dM_s\right] = 0.$$

So,

$$\mathbb{E}[f(t \wedge T_n, X_{t \wedge T_n})]$$
$$= \mathbb{E}\left[\int_0^{t \wedge T_n} D_x f(s, X_{s-}) \, dA_s\right]$$
(29)
$$+ \frac{1}{2}\mathbb{E}\left[\int_0^{t \wedge T_n} \theta(s, X_s) \, d[X]_s^c\right] + \mathbb{E}\left[\int_0^{t \wedge T_n} D_t f(s, X_s) \, ds\right]$$
$$+ \mathbb{E}\left[\sum_{s \leq t \wedge T_n} \int_{X_{s-}}^{X_s} (X_s - x)\theta(s,x) \, dx\right].$$

Letting $n$ go to infinity, monotone convergence implies convergence of the second and fourth terms on the right-hand side and dominated convergence implies convergence of the first and third terms. Also, uniform integrability of $X_{t \wedge T_n} = M_{t \wedge T_n} + A_{t \wedge T_n}$ over $n \in \mathbb{N}$ shows that the term on the left-hand side will also converge.

Taking $t$ large enough so that the support of $\theta$ is contained in $[0,t] \times \mathbb{R}$, $f(t, X_t) = 0$ and taking the limit as $n \to \infty$ in (29) gives

$$0 = \mathbb{E}\left[\int D_x f(s, X_{s-}) \, dA_s\right] + \frac{1}{2}\mathbb{E}\left[\int \theta(s, X_s) \, d[X]_s^c\right]$$
$$+ \mathbb{E}\left[\int D_t f(s, X_s) \, ds\right] + \mathbb{E}\left[\sum_{s>0} \int_{X_{s-}}^{X_s} (X_s - x)\theta(s,x) \, dx\right].$$

The result now follows by substituting in

$$\mathbb{E}\left[\int D_x f(s, X_{s-})\, dA_s\right] = \mathbb{E}\left[\int\int_{-\infty}^{X_{s-}} \theta(s,y)\, dy\, dA_s\right]$$

and by using integration by parts

$$\int \mathbb{E}[D_t f(s, X_s)]\, ds = \int \mathbb{E}\left[\int D_t \theta(s,x)(X_s - x)_+\, dx\right] ds$$
$$= \int\int D_t \theta(t,x) C(t,x)\, dt\, dx$$
$$= -\int\int \theta(t,x)\, d_t C(t,x)\, dx. \qquad \square$$

The following simple consequence of Lemma 3.3 will be used to show that (10) is satisfied.

COROLLARY 3.4. *Let $X$ be a càdlàg adapted process which decomposes as $X = M + A$ for a martingale $M$ and càdlàg integrable process $A$ with integrable variation over finite time intervals. Define $C \in \mathcal{D}$ by $C(t,x) \equiv \mathbb{E}[(X_t - x)_+]$.*

*If $f(t,x)$ is locally Lipschitz continuous in $x$ then,*

$$\mathbb{E}\left[\int 1_{\{(t, X_t) \notin \operatorname{diff}(f)\}}\, d[X]_t^c\right] = 2\int\int 1_{\{(t,x) \notin \operatorname{diff}(f)\}} |d_t C(t,x)|\, dx.$$

PROOF. First, choose any nonnegative bounded and measurable $\theta \colon \mathbb{R}_+ \times \mathbb{R} \to \mathbb{R}$ with bounded support.

We use a result of Lebesgue which states that any locally Lipschitz continuous function on the reals is differentiable almost everywhere ([13], Theorem 3.2), giving

$$\int_{-\infty}^{X_{t-}} 1_{\{(t,y) \notin \operatorname{diff}(f)\}} \theta(t,y)\, dy = \int_{X_{t-}}^{X_t} 1_{\{(t,x) \notin \operatorname{diff}(f)\}} (X_t - x)\theta(t,x)\, dx = 0.$$

So, replacing $\theta(t,x)$ by $1_{\{(t,x) \notin \operatorname{diff}(f)\}}\theta(t,x)$ in (28) gives

$$\mathbb{E}\left[\int 1_{\{(t, X_t) \notin \operatorname{diff}(f)\}} \theta(t, X_t)\, d[X]_t^c\right] = 2\int\int 1_{\{(t,x) \notin \operatorname{diff}(f)\}} \theta(t,x)\, d_t C(t,x)\, dx.$$

Letting $\theta$ increase to 1 gives the result. $\square$

We can now complete the proof of Theorem 1.2, which makes use of Lebesgue's result that locally Lipschitz continuous functions of the reals are differentiable almost everywhere ([13], Theorem 3.2).

PROOF OF THEOREM 1.2. As in the proof of Lemma 3.2, pre-localization can be used to reduce to the case where $X$ decomposes as the sum of a martingale and a càdlàg integrable process with integrable variation. Define $C \in \mathcal{D}$ by $C(t,x) = \mathbb{E}[(X_t - x)_+]$.

Let $\text{diff}(f)$ be the set of $x \in \mathbb{R}$ at which $f$ is differentiable. Also, choose any $t > 0$ and set $\theta(x) = \int_0^t |d_s C(s,x)|$. Corollary 3.4 gives

$$\mathbb{E}\left[\int_0^t 1_{\{X_s \notin \text{diff}(f)\}} \, d[X]_s^c\right] = 2 \int 1_{\{x \notin \text{diff}(f)\}} \theta(x) \, dx.$$

However, as $f$ is locally Lipschitz continuous, Lebesgue's theorem tells us that $f$ is differentiable almost everywhere, and the right-hand side of the above equality is 0. So, (10) is satisfied and Lemma 3.2 gives (11). The result now follows from Theorem 2.1. $\square$

The proof of Theorem 1.4 also follows easily.

PROOF OF THEOREM 1.4. As above, we may restrict to the case where $X$ is a sum of a martingale and a càdlàg integrable process with integrable variation. Then Corollary 3.4 and Lemma A.3 give

$$\mathbb{E}\left[\int 1_{\{(t,X_t) \notin \text{diff}(f)\}} \, d[X]_t^c\right] = 2 \iint 1_{\{(t,x) \notin \text{diff}(f)\}} |d_t C(t,x)| \, dx = 0.$$

Therefore, (10) is satisfied and Lemma 3.2 gives (11), so Theorem 2.1 gives the result. $\square$

## APPENDIX: "ALMOST EVERYWHERE" DIFFERENTIABILITY

In this appendix, we show that functions in $\mathcal{D}$ are differentiable in the "almost everywhere" sense required by the proof of Theorem 1.4. See Lemma A.3 below for the statement of the result.

For every $a \in \mathbb{R} \setminus \{0\}$ we use $\nabla_a$ to represent the finite difference operator

$$\nabla_a f(t,x) \equiv (f(t, x+a) - f(t,x))/a.$$

Also, for $f \in \mathcal{D}_0$, the left limit in $t$ is denoted by

$$f^-(t,x) \equiv \begin{cases} \lim_{s \uparrow\uparrow t} f(s,x), & \text{if } t > 0, \\ f(0,x), & \text{if } t = 0. \end{cases}$$

Then we have the following integration by parts formula.

LEMMA A.1. *Suppose that one of $f, g \in \mathcal{D}_0$ has compact support in $(0, \infty) \times \mathbb{R}$. Then*

(30) $$\iint \nabla_a f \, d_t g \, dx = \iint \nabla_{-a} g^- \, d_t f \, dx.$$

PROOF. Choosing any $a \in \mathbb{R}$, integration by parts and the condition that $fg$ has compact support in $(0,\infty) \times \mathbb{R}$, gives

$$\int f(t, x+a)\, d_t g(t,x) + \int g^-(t,x)\, d_t f(t, x+a) = 0.$$

Then integrate with respect to $x$,

$$\iint f(t, x+a)\, d_t g(t,x)\, dx + \iint g^-(t, x-a)\, d_t f(t,x)\, dx = 0.$$

The result follows by subtracting this equation from itself with $a$ replaced by 0 and dividing by $a$. □

Letting $\hat{\nabla}_a f$ be the difference of the left and right finite differences

$$\hat{\nabla}_a f \equiv \nabla_a f - \nabla_{-a} f,$$

then the previous result can be used to prove the following limit.

LEMMA A.2. *Let $f, g \in \mathcal{D}_0$ and $\theta : \mathbb{R}_+ \times \mathbb{R} \to \mathbb{R}$ be measurable and bounded with bounded support. Then*

$$\iint \theta \hat{\nabla}_a f\, d_t g\, dx + \iint \theta \hat{\nabla}_a g\, d_t f\, dx \to 0 \tag{31}$$

*as $a \to 0$.*

PROOF. First, if we suppose that $f$ has compact support in $(0, \infty) \times \mathbb{R}$, then we can take the difference of (30) with itself, with $a$ replaced by $-a$ to get

$$\iint \hat{\nabla}_a f\, d_t g\, dx + \iint \hat{\nabla}_a g^-\, d_t f\, dx = 0. \tag{32}$$

Now for general $f, g \in \mathcal{D}_0$, choose any continuously differentiable $\theta : \mathbb{R}_+ \times \mathbb{R} \to \mathbb{R}$ with compact support in $(0, \infty) \times \mathbb{R}$. Replacing $f$ by $\theta f$ in (32),

$$\iint \hat{\nabla}_a(\theta f)\, d_t g\, dx + \iint \theta \hat{\nabla}_a g^-\, d_t f\, dx + \iint (\hat{\nabla}_a g^-) \frac{\partial \theta}{\partial t} f\, dt\, dx = 0. \tag{33}$$

For any function $h(t,x)$ which is locally Lipschitz continuous in $x$, we can make use of the identity $\nabla_{-a} h(t,x) = \nabla_a h(t, x-a)$ to get

$$\int u(x) \hat{\nabla}_a h(t,x)\, dx$$
$$= \int (u(x) \nabla_a h(t,x) - u(x) \nabla_a h(t, x-a))\, dx$$
$$= \int (u(x) - u(x+a)) \nabla_a h(t,x)\, dx \to 0$$

FUNCTIONS OF SEMIMARTINGALES 25as $a \to 0$, whenever $u$ is continuous with compact support. Combining this with dominated convergence for the following integrals gives

$$\iint (\hat{\nabla}_a g^-) \frac{\partial \theta}{\partial t} f \, dt \, dx = \int \left( \int (\hat{\nabla}_a g^-) \frac{\partial \theta}{\partial t} f \, dx \right) dt \to 0,$$

$$\iint \theta(\hat{\nabla}_a(g^- - g)) \, d_t f \, dx = \sum_{t>0} \int \theta \hat{\nabla}_a (g^- - g)(f - f^-) \, dx \to 0$$

as $a \to 0$. In the second of these limits, the fact that that there are only countably many times at which $g^- \neq g$ has been used to write the integral as an infinite sum. Combining these limits with (33),

(34) $$\iint \hat{\nabla}_a(\theta f) \, d_t g \, dx + \iint \theta \hat{\nabla}_a g \, d_t f \, dx \to 0$$

as $a \to 0$. Now, the limit

$$\nabla_a(\theta f) - \theta \nabla_a f = f \nabla_a \theta + a \nabla_a \theta \nabla_a f \to f \frac{\partial \theta}{\partial x}$$

as $a \to 0$ implies that $\hat{\nabla}_a(\theta f) - \theta \hat{\nabla}_a f \to 0$. Applying this with dominated convergence to the first integral in (34) gives (31). The result for arbitrary $\theta$ follows from the monotone class lemma. $\square$

Finally for this section, Lemma A.2 is used to show that every $f \in \mathcal{D}$ is differentiable with respect to $x$ in an "almost everywhere" sense. It is not clear if this result will generalize to arbitrary $f \in \mathcal{D}_0$ which would imply that Theorem 1.4 holds for all $f$ in $\mathcal{D}_0$.

LEMMA A.3. *For any $f \in \mathcal{D}$ and $g \in \mathcal{D}_0$,*

$$\iint 1_{\{(t,x) \notin \mathrm{diff}(f)\}} |d_t g(t,x)| \, dx = 0.$$

PROOF. Write $D_x^- f(t,x)$ and $D_x^+ f(t,x)$ for the left and right derivatives, respectively, of $f(t,x)$ by $x$. Then $f(t,x)$ is differentiable with respect to $x$ at those points where $D_x^- f = D_x^+ f$.

The definition of $\hat{\nabla}_a$ gives $\hat{\nabla}_a f \to D_x^+ f - D_x^- f$ as $a \downarrow 0$. So (31) with $g$ replaced by $f$ gives

$$\iint \theta(D_x^+ f - D_x^- f) \, d_t f \, dx = 0.$$

As this is true for every bounded and measurable $\theta$ with bounded support,

(35) $$\iint 1_{\{(t,x) \notin \mathrm{diff}(f)\}} |d_t f| \, dx = 0.$$

Similarly, (31) gives

$$\text{(36)} \qquad \iint \theta(D_x^+ f - D_x^- f)\, d_t g\, dx = -\lim_{a \downarrow 0} \iint \theta \hat{\nabla}_a g\, d_t f\, dx.$$

Letting $K(t,x)$ be the locally bounded function $\limsup_{a \downarrow 0} |\hat{\nabla}_a g(t,x)|$, applying dominated convergence to the right-hand side of (36) gives

$$\left| \iint \theta(D_x^+ f - D_x^- f)\, d_t g\, dx \right| \le \iint |\theta| K |d_t f|\, dx.$$

Replace $\theta$ by $1_{\{(t,x) \notin \text{diff}(f)\}} \theta$ in this inequality and apply (35),

$$\left| \iint \theta(D_x^+ f - D_x^- f)\, d_t g\, dx \right| \le \iint 1_{\{(t,x) \notin \text{diff}(f)\}} |\theta| K |d_t f|\, dx = 0.$$

As this is true for every measurable and bounded $\theta$ with bounded support, the result follows. □

6 St Peter's Street
London, N1 8JG
United Kingdom
E-mail: george.lowther@blueyonder.co.uk